\documentclass[oneside,12pt]{amsart}
\usepackage[all]{xy}

\usepackage{amscd,amsmath,latexsym, amssymb}
\usepackage[mathcal]{euscript}

\newtheorem{thm}{Theorem}[section]

\newtheorem{lem}[thm]{Lemma}

\theoremstyle{definition}
\newtheorem{defin}[thm]{Definition}
\theoremstyle{definition}

\newtheorem{remark}[thm]{Remark}
\theoremstyle{remark}

\begin{document}
\title[Cup-products in generalized moment-angle complexes] {Cup-products in generalized moment-angle complexes}

\author[A.~Bahri]{A.~Bahri}
\address{Department of Mathematics,
Rider University, Lawrenceville, NJ 08648, U.S.A.}
\email{bahri@rider.edu}

\author[M.~Bendersky]{M.~Bendersky}
\address{Department of Mathematics
CUNY,  East 695 Park Avenue New York, NY 10065, U.S.A.}
\email{mbenders@xena.hunter.cuny.edu}

\author[F.~R.~Cohen]{F.~R.~Cohen}
\address{Department of Mathematics,
University of Rochester, Rochester, NY 14625, U.S.A.}
\email{cohf@math.rochester.edu}

\author[S.~Gitler]{S.~Gitler}
\address{Department of Mathematics,
Cinvestav, San Pedro Zacatenco, Mexico, D.F. CP 07360 Apartado
Postal 14-740, Mexico} \email{sgitler@math.cinvestav.mx}

\subjclass{Primary:  13F55, 14F45, 32S22, 52C35, 55U10 Secondary:
16E05, 57R19}

\keywords{arrangements, cohomology ring, moment-angle complex,
simplicial complex, simplicial sets, stable splittings,
Stanley-Reisner ring, suspensions, toric varieties}

\begin{abstract}
Given a family of based CW-pairs
$(\underline{X},\underline{A})=\{(X_i,A_i)\}^m_{i=1}$ together with
an abstract simplicial complex $K$ with $m$ vertices, there is an
associated based CW-complex $Z(K;(\underline{X},\underline{A}))$
known as a generalized moment-angle complex \cite{bbcg}.

The decomposition theorem of \cite{bbcg}, \cite{bbcg2} splits the
suspension of $Z(K; (\underline{X}, \underline{A}))$ into  a bouquet
of spaces determined by the full sub-complexes of $K$. That
decomposition theorem is used here to describe the ring structure
for the cohomology of Z(K; (\underline{X}, \underline{A})). Explicit
computations are made for families of suspension pairs and for the
cases where $X_i$ is the cone on $A_i$. These results complement and
generalize those of Davis-Januszkiewicz \cite{davis.jan}, Franz,
\cite{franz} and \cite{franz2}, Hochster \cite{hochster} as well as
Panov  \cite{panov} and Baskakov-Buchstaber-Panov, \cite{bbp}. Under
conditions stated below, these theorems also apply for generalized
cohomology theories.
\end{abstract}

\maketitle

\section{\bf Introduction, definitions, and main results}\label{Introduction}

This paper is a study of the cup-product structure for the
cohomology ring of a generalized moment-angle complex. The new
result here is that the structure of the cohomology ring is given in
terms of a geometric decomposition arising after one suspension of
the generalized moment-angle complex \cite{bbcg, bbcg2}.

This cup-product structure was studied for special cases in
\cite{davis.jan, bbp, buchstaber.panov, hochster,
panov,franz,franz2} with important basic cases first given by
Davis-Januszkiewicz \cite{davis.jan}, Franz \cite{franz, franz2},
and Buchstaber-Panov \cite{buchstaber.panov,panov}. A few details
concerning historical developments are listed in \cite{bbcg, bbcg2}.
The methods here give a determination of the cohomology ring
structure for many new generalized moment-angle complexes as well as
retrieve many known results.

A generalized moment-angle complex is a union of cartesian products
of based CW-complexes \cite{bbcg, bbcg2}. Generalized moment-angle
complexes satisfy a stable decomposition extending a classical
decomposition for suspensions of product spaces. This decomposition
informs on the cup product structure by providing information about
the diagonal map, after stabilization; this suffices to give the
cup-product structure. The notation used in this article is adopted
from \cite{bbcg, bbcg2}.

Consider a product of based CW-complexes
$Y^{[m]}=\prod\limits^m_{i=1}Y_i$. Let $I\subseteq[m]$ be an
increasing subsequence of $[m]=(1,2,\ldots,m)$. If
$I=(i_1,\ldots,i_k)$, then $\widehat{Y}^I$ denotes the smash product
$Y_{i_1}\wedge\cdots\wedge Y_{i_k}$ the quotient space of
$Y^I=Y_{i_1} \times\cdots \times Y_{i_k}$ by the subspace given by
the fat wedge $FW(Y^I)=\{(y_{i_1},\ldots , y_{i_k})\in Y^I\quad
|\quad y_{i_j}=$ base-point of $Y_{i_j}$ for at least one $i_j$)\}.

Recall that a generalized moment-angle complex is a functor of two
variables:
\begin{enumerate}
  \item An abstract simplicial complex $K$ with $m$ vertices identified with the
sequence $(1,\ldots,m)$. Then the simplices of $K$ are identified
with increasing subsequences $\sigma=(i_1,\ldots,i_t)$. The
dimension of $\sigma$ is $t-1$. The defining property of $K$ is that
if $\tau\subset\sigma$ is a subsequence of $\sigma$ then $\tau\in
K$. The empty set $\varnothing$ belongs to $K$.
\item A family $(\underline{X},\underline{A})=\{(X_i,A_i,x_i)\}_{i\in[m]}$ of
connected, based CW-pairs $(X_i,A_i,x_i)$.
\end{enumerate}

The morphisms in (1) are embeddings of simplicial complexes and the
morphisms in (2) are maps of based connected pairs
$(\underline{X},\underline{A},\underline{x})\to
(\underline{Y},\underline{B}, \underline{y})$. The simplicial
complex $K$ has a family of full sub-complexes $K_I$ defined for
every subsequence $I$ of $[m]$,
$$K_I=\{\sigma\cap I|\sigma\in K\}$$
$$K_I\quad{\rm has}\quad l(I),\,\,\hbox{\rm length of }\,\, I,\quad
{\rm vertices},$$ and associated to $K_I$ a family of spaces
$(\underline{X},\underline{A})_I=\{(X_i,A_i)\}_{i\in I}$. If
$I=[m]$, then $K_I=K$.

Next, define the functors $Z(K;(\underline{X},\underline{A}))$ and
$\widehat{Z}(K;(\underline{X},\underline{A}))$ as in
\cite{bbcg,bbcg2} in the following way. For every $\sigma\in K$,
define
$$
D(\sigma) =\prod^m_{i=1}Y_i,\quad {\rm where}\quad Y_i
=\left\{\begin{array}{lcl}
X_i &{\rm if} & i \in \sigma\\
A_i &{\rm if} & i \in [m]-\sigma.
\end{array}\right.$$ with $D(\varnothing) = A_1 \times \ldots \times A_m$.
As in \cite{bbcg}, the space $Z(K;(\underline{X},\underline{A}))$ is
defined as $Z(K;(\underline{X},\underline{A}))=\bigcup_{\sigma \in
K} D(\sigma)= \operatorname{colim} D(\sigma).$

In what follows, it is useful to define variations for a fixed,
ambient $I$ where the analogue of $D(\sigma)$ is replaced as
follows.
\begin{defin}\label{defin:moment.angle.complex.variation}
For fixed $I=(i_1,\ldots,i_k)$, and every $\sigma\in K$, define
\begin{equation}
\label{eq: Y.I 1x1111}
    \begin{array}{l}
    Y^I(\sigma\cap I)=Y_{i_1} \times\cdots \times Y_{i_k},
        \end{array}
\end{equation}and

\begin{equation}
\label{eq: Y.smash.I 1x1111}
    \begin{array}{l}
    \widehat{Y}^I(\sigma\cap I)=Y_{i_1}\wedge\cdots\wedge Y_{i_k}
        \end{array}
\end{equation} where

\begin{equation}
\begin{array}{l}
Y_j=\left\{\begin{array}{l}
X_j\quad{\rm if}\quad j \in\sigma \cap I\\
A_j\quad{\rm if}\quad j \in I-\sigma\cap I.
\end{array}\right.\\[-.5cm]
\end{array}
\end{equation}

Furthermore,
\begin{eqnarray*}
Y^I(\Phi)&=&A^I=A_{i_1}\times\cdots\times A_{i_k}\\
\widehat{Y}^I(\Phi)&=&\widehat{A}^I=A_{i_1}\wedge\cdots\wedge
A_{i_k}.
\end{eqnarray*}
Then the generalized moment-angle complexes are
$$Z(K_I;(\underline{X},\underline{A})_I)=\bigcup_{\sigma \in K} Y^I(\sigma\cap I)$$
and
$$\widehat{Z}(K_I;(\underline{X},\underline{A})_I)=\bigcup\limits_{\sigma
\in K}\widehat{Y}^I(\sigma\cap I).$$ {\bf Note:} The notation
$Z(K_I;(\underline{X_I},\underline{A_I}))$ was used in \cite{bbcg}
for $Z(K_I;(\underline{X},\underline{A})_I)$.

To simplify notation below, the following notation
$$ Z(K),\ \widehat{Z}(K), \  Z(K_I) \;{\rm  and}\; \widehat{Z}(K_I) $$ is used to denote
$Z(K;(\underline{X},\underline{A})), \
\widehat{Z}(K;(\underline{X},\underline{A})), \
Z(K_I;(\underline{X},\underline{A})_I),\ \;{\rm  and}\;
\widehat{Z}(K_I;(\underline{X},\underline{A})_I)$ respectively.

\end{defin}

The results of \cite{bbcg} stated next are the main ingredients used
here to analyze the cup-product structure for the generalized
moment-angle complex.
\begin{thm} \label{thm:bbcg}
Let $K$ be an abstract simplicial complex with $m$ vertices. Assume
that $$(\underline{X},\underline{A}) =\{(X_i, A_i,x_i)\}^m_{i=1}$$
are pointed triples of CW-complexes for all $i$. Then there is a
natural, pointed homotopy equivalence
$$H: \Sigma(Z(K;(\underline{X},\underline{A})))\to \Sigma(\bigvee_{I \subseteq [m]}
\widehat{Z}(K_I;(\underline{X},\underline{A})_I)).$$
\end{thm}

Cup-products in the cohomology of any space $W$ are induced by the
diagonal map $$W \to W \times W.$$ The main direction of this paper
is an analysis of the behavior of the diagonal map for the
generalized moment-angle complex and the properties of the diagonal
map which are preserved by the stable decomposition of Theorem
\ref{thm:bbcg} above \cite{bbcg}.

Let $$\Delta_I:Y^I\to Y^I\wedge Y^I$$ denote the reduced diagonal of
$Y^I$ and let $$\widehat{\Delta}_I:
\widehat{Y}^I\to\widehat{Y}^I\wedge\widehat{Y}^I$$ denote the
reduced diagonal of $\widehat{Y}^I$. In this paper {\it partial
diagonals} are defined below
$$\widehat{\Delta}^{J,L}_I:\widehat{Y}^I\to \widehat{Y}^J\wedge\widehat{Y}^L,$$
and by restriction $$\widehat{\Delta}^{J,L}_I:\widehat{Z}(K_I)\to
\widehat{Z}(K_J)\wedge\widehat{Z}(K_L)$$ where $J\cup L=I$. If
$I=J=L$, these maps coincide with the {\it reduced} diagonal
$\widehat{\Delta}_I$. Furthermore, if
$\widehat{\Pi}_I:Y^{[m]}\to\widehat{Y}^I$ is the projection, there
are commutative diagrams of $CW$-complexes and based continuous maps

\begin{equation}
\label{eq:1.1} \xymatrix{ Y^{[m]} \ar[r]^{\kern-.6cm\Delta_{[m]}}
\ar[d]_{\widehat{\Pi}_I} &Y^{[m]}\wedge Y^{[m]}
\ar[d]^{\widehat{\Pi}_J\wedge\widehat{\Pi}_L}\\
\widehat{Y}^I\ar[r]^{\kern-.6cm\widehat{\Delta}^{J,L}_I}&\widehat{Y}^J\wedge\widehat{Y}^L}
\end{equation} and by restriction to $Z(K) \subset X^{[m]}$

\begin{equation}
\label{eq:1.2} \xymatrix{ Z(K) \ar[r]^{\kern-.7cm\Delta_K}
\ar[d]_{\widehat{\Pi}_I} &Z(K)\wedge Z(K)
\ar[d]^{\widehat{\Pi}_J\wedge\widehat{\Pi}_L}\\
\widehat{Z}(K_I)\ar[r]^{\kern-.7cm\widehat{\Delta}^{J,L}_I}&
\widehat{Z}(K_J)\wedge\widehat{Z}(K_L).}
\end{equation}

A definition is given next.
\begin{defin}\label{defin:star.product}
Assume given a family of based CW-pairs
$(\underline{X},\underline{A})=\{(X_i,A_i)\}^m_{i=1}$. Given
cohomology classes $u\in H^p(Z(K_J)), v\in H^q(Z(K_L))$, define
$$u*v=(\widehat{\Delta}^{J,L}_I)^*(u\otimes v)\quad\hbox{\rm thus}\quad u*v\in H^{p+q}(\widehat{Z}(K_I)).$$
The element $u*v\in H^{p+q}(\widehat{Z}(K_I))$ is called the
$*$-product. Commutativity of diagram (\ref{eq:1.2}) gives
\begin{equation}
\label{eq:1.3} \widehat{\Pi}^*_I(u*v)=
\widehat{\Pi}^*_J(u)\smallsmile \widehat{\Pi}^*_L(v)
\end{equation} where $\smallsmile$ denotes the cup-product for the CW-complex $Z(K)$.

Let $$\mathcal H^q(K;(\underline{X},\underline{A}))=
\bigoplus_{I\subseteq m} H^q(\widehat{Z}(K_I))$$ with
$$\mathcal H^*(K;(\underline{X},\underline{A}))=
\bigoplus_{I\subseteq m}H^*(\widehat{Z}(K_I)).$$ Define a map
$$ \eta:\mathcal H^*(K;(\underline{X},\underline{A}))\to H^*({Z}
(K;(\underline{X},\underline{A}))$$%
where $\eta$ restricted to $H^*(\widehat{Z}(K_I))$ is
$\widehat{\Pi}^*_I$.
\end{defin}

By the decomposition given in Theorem $2.8$ of \cite{bbcg},
$\eta=\bigoplus\limits_{I\leq [m]}\Pi^*_I$ is an additive
isomorphism. The $*$-product gives $\mathcal H
^*(K;(\underline{X},\underline{A}))$ the structure of an algebra, a
fact which is checked in Section \ref{sec:The partial diagonal in
smash moment angle} where the next result is proven.

\begin{thm}\label{thm:1.4}
Let $K$ be an abstract simplicial complex with $m$ vertices. Assume
that $(\underline{X},\underline{A})=\{(X_i,A_i, x_i)\}^m_{i=1}$ is a
family of based CW-pairs. Then
$$\eta:\mathcal H^*(K; (\underline{X},\underline{A}))\to H^*(Z(K;
(\underline{X},\underline{A}))$$ is a ring isomorphism.
\end{thm}

\begin{defin}\label{defin:suspension.pairs}
Assume that $(\underline{X},\underline{A})=\{(X_i,A_i,
x_i)\}^m_{i=1}$ is a family of based CW-pairs. The pair
$(\underline{X},\underline{A})$ is a {\it suspension pair} if
$(X_i,A_i) = (\Sigma(U_i), \Sigma(V_i))$ for each $i$ with each
inclusion $A_i \subset X_i$  given as the suspension of a map
$f_i:V_i \to U_i$.
\end{defin}

If the pair $(\underline{X},\underline{A})$ is a suspension pair,
then the reduced diagonal $$\Delta_i:Y_i\to Y_i\wedge Y_i$$ for $Y_i
= X_i \ \mbox{or} \ A_i$ is null-homotopic. This fact will be used
below to prove the next result.
\begin{thm}\label{thm:1.5}
Let $K$ be an abstract simplicial complex with $m$ vertices. Assume
that $(\underline{X},\underline{A})=\{(X_i,A_i,x_i)\}^m_{i=1}$ is a
family of based CW-pairs.  If $(\underline{X},\underline{A})
=\{(\Sigma(U_i),\Sigma( V_i))\}_{i\in[m]}$ is a suspension pair, and
$J \cap L \neq \varnothing$, then
$$\widehat{\Delta}^{J,L}_I:\widehat{Z}(K_I)\to
\widehat{Z}(K_J)\wedge\widehat{Z}(K_L)$$ is null-homotopic, and thus
$u*v=0$ for classes $u\in H^p\widehat{Z}(K_J)$ and $v\in
H^q\widehat{Z}(K_L)$.
\end{thm}

\begin{defin}\label{defin:stably.wedge.equivalent}
Define two CW-complexes $X$ and $Y$ to be {\it stably wedge
equivalent}, if (i) $X$ is stably equivalent to $X_1\vee\cdots\vee
X_k$, (ii) $Y$ is stably equivalent to $Y_1\vee\cdots\vee Y_k$, and
(iii) $X_i$ is stably equivalent to $Y_i$ for $i=1,\ldots,k$. In
particular, if $X$ and $Y$ are stably wedge equivalent, then they
are stably homotopy equivalent. Let $T=(t_1,\ldots,t_m)$ be a
sequence of positive integers, then define
$$(\underline{\Sigma^TX}, \underline{\Sigma^TA})=\left\{(\Sigma^{t_i}X_i,\Sigma^{t_i}A_i)\right\}^m_{i=1}$$
where $\Sigma^{t_i}X_i$ is the $t_i$-th iterated suspension of
$X_i$. Let $T=(t_1,\ldots,t_m)$, and $T'=(t'_1,\ldots,t'_m)$ denote
two sequences of strictly positive integers of length $m$, and
define $T\equiv T'$ mod $2$ if $t_i\equiv t'_i$ mod $2$ for all $i$.
\end{defin}

\begin{thm}\label{thm:1.6} Let $K$ be an abstract simplicial complex with $m$ vertices. Assume
that $(\underline{X},\underline{A})=\{(X_i,A_i,x_i)\}^m_{i=1}$ is a
family of based CW-pairs and that $T=(t_1,\ldots,t_m)$, and
$T'=(t'_1,\ldots,t'_m)$ are two sequences of strictly positive
integers of length $m$. Then the following hold.
\begin{enumerate}
\item $Z(K;(\underline{\Sigma^TX}, \underline{\Sigma^TA}))$ and
$Z(K;(\underline{\Sigma^{T'}X},\underline{\Sigma^{T'}A}))$ are
stably wedge equivalent.

\item If $T\equiv T'$ mod $2$, then $Z(K;(\underline{\Sigma^TX}, \underline{\Sigma^TA}))$
and $Z(K;(\underline{\Sigma^{T'}X},\underline{\Sigma^{T'}A}))$ have
isomorphic cohomology rings regarded as ungraded rings.
\end{enumerate}
\end{thm}

\begin{thm}\label{thm:cup.products.for.CX.X}
Let $K$ be an abstract simplicial complex with $m$ vertices. Assume
that $(\underline{CX},\underline{X})=\{(CX_i,X_i,x_i)\}^m_{i=1}$ is
a family of based CW-pairs such that the finite product \
$$(X_{1} \times \cdots  \times X_{m})\times (Z(K_{I_1};(D^1,S^0))
\times \cdots \times Z(K_{I_t};(D^1,S^0)))$$ for all $I_j \subseteq
[m]$ satisfies the strong form of the K\"unneth theorem. Then the
cup-product structure for the cohomology algebra
$H^*(Z(K;(\underline{CX},\underline{X})))$ is a functor of the
cohomology algebras of $X_i  \mbox{ for all} \ i$, and
$Z(K_I;(D^1,S^0))$ for all $I$.
\end{thm}

The analogue of Theorem  \ref{thm:cup.products.for.CX.X} in the case
of $H^*(Z(K;(\underline{X},\underline{A})))$ for which $A$ is
contractible is given as Theorem $2.35$ of \cite{bbcg}. In this
case, the result depends only on the structure of the cohomology
algebra of $H^*(X_i)$, $1 \leq i \leq m$.

The cohomology ring of $Z(K;(D^2,S^1))$ was studied by Panov
\cite{panov} as well as others \cite{bbp,buchstaber.panov,franz2}.
Panov proved that the cohomology ring was the Tor-algebra of the
face ring of $K$ as studied by Hochster \cite{hochster}. Theorem
\ref{thm:cup.products.for.CX.X} gives information about that ring.

\begin{remark}\label{remark:generalized.cohomology.theory}
Theorems \ref{thm:1.4}, \ref{thm:1.5}, and \ref{thm:1.6} above apply
more generally to any cohomology theory. On the other hand, Theorem
\ref{thm:cup.products.for.CX.X} applies in the case of any
cohomology theory which satisfies additional restrictions formalized
as follows.
\end{remark}

\begin{defin}\label{defin:good.spaces.withy respect.to.a.cohomology.theory}
A family of based CW-pairs
$(\underline{X},\underline{A})=\{(X_i,A_i,x_i)\}^m_{i=1}$ together
with a finite simpicial complex $K$ with $m$ vertices is said to be
{\bf proper with respect to a cohomology theory $E^*(-)$} provided
the strong form of the K\"unneth theorem is satisfied for any finite
smash product of the spaces given by
\begin{enumerate}
  \item $X_i$,
  \item $A_i$, and
  \item $|K_I|$ for all $I \subseteq [m]$.
\end{enumerate}
\end{defin}

The analogue of Theorem \ref{thm:cup.products.for.CX.X} for any
cohomology theory follows next.
\begin{thm}\label{thm:E.cup.products.for.CX.X}
Let $(\underline{CX},\underline{X})=\{(CX_i,X_i)\}^m_{i=1}$ be a
family of based CW-pairs, and $K$ be an abstract simplicial complex
with $m$ vertices which is proper with respect to a cohomology
theory $E^*(-)$. Then the cup-product structure for the cohomology
algebra $E^*(Z(K;(\underline{CX},\underline{X})))$ is a functor of
\begin{enumerate}
\item the cohomology algebras $E^*(X_i)$, $1 \leq i \leq m$, and
\item the cohomology algebras $E^*(Z(K_I;(D^1,S^0)))$ for all $I\subseteq K$.
\end{enumerate}

Furthermore, there are natural isomorphisms of multiplicatively
closed sub-modules of
$$E^*(\Sigma|K_I|)) \otimes E^*(X_1) \otimes \cdots \otimes E^*(X_m)$$ given by
$$\tilde{E}^*(\widehat{Z}(K_I;(\underline{CX},\underline{X})_I)) \to
\tilde{E}^*(\Sigma|K_I|)) \otimes \tilde{E}^*(\widehat{X}^I),$$ and
$$\tilde{E}^*(\Sigma|K_I|)) \otimes \tilde{E}^*(\widehat{X}^I) \to \tilde{E}^*(\Sigma|K_I|)) \otimes \tilde{E}^*(X_{i_1}) \otimes
\tilde{E}^*(X_{i_2}) \otimes \cdots \otimes \tilde{E}^*(X_{i_k})$$
for $I= (i_1, \ldots, i_k) \subseteq [m]$.

In addition, there are natural isomorphisms
$$ \underset{I \subseteq [m]}{\bigoplus }\tilde{E}^*(
\widehat{Z}(K_I;(\underline{CX},\underline{X})_I)) \to
\widetilde{E}^*(Z(K;(\underline{CX},\underline{X})))$$ where
$\tilde{E}^*( \widehat{Z}(K_I;(\underline{CX},\underline{X})_I))$ is
isomorphic to a multiplicatively closed sub-module of
$$\widetilde{E}^*(Z(K;(\underline{CX},\underline{X}))).$$
\end{thm}

The definition of {\it partial diagonals} is the subject of section
\ref{sec:The partial diagonal in product spaces}. Information about
partial diagonals is extended to smash product generalized
moment-angle complexes in section \ref{sec:The partial diagonal in
smash moment angle}. Sections \ref{sec:Proof of Theorem (1.4)},
\ref{sec:Proof of Theorem (1.5)}, and \ref{sec:Proof of Theorem
(1.6)} give the proofs of Theorems \ref{thm:1.4}, \ref{thm:1.5}, and
\ref{thm:1.6} respectively. 

Partial diagonals are given in section \ref{sec:CX.X}
for the smash moment-angle complexes
$\widehat{Z}(K;(\underline{CX},\underline{X})$ where $CX$ is the
cone on the CW-complex $X$. In this case, the cohomology algebra of
$Z(K;(CX,X))$ is shown to be a functor, with mild restrictions, of
the cohomology algebra for $X$ and for $Z(K;(D^1,S^0))$. This result
is a counter-point to Theorem $2.35$ of \cite{bbcg} where an
analogous result is proven for $Z(K;(X,*))$.

\tableofcontents

\section{\bf The partial diagonal in product spaces}\label{sec:The partial diagonal in product spaces}

Let $Y^{[m]}=Y_1\times\cdots\times Y_m$ and
$\widehat{Y}^I=Y_{i_j}\wedge\cdots\wedge Y_{i_k}$ for
$I=(i_1,\ldots,i_k)\subseteq [m]$. There are natural projection maps
$\widehat{\Pi}_I:Y^{[m]}\to \widehat{Y}^I$ obtained as the
composition
$$Y^{[m]}\stackrel{\Pi_I}{\to}Y^I\stackrel{\rho_I}{\to}\widehat{Y}^I$$ where $\Pi_I$ is the projection map and $\rho_I$ is the quotient map.

Let $$\widehat{\Delta}_I:\widehat{Y}^I\to
\widehat{Y}^I\wedge\widehat{Y}^I$$ be the reduced diagonal map of
$\widehat{Y}^I$, and define $$C_I=\{(J,L)\  | \ J,L\subseteq
I\quad{\rm and}\quad J\cup L=I\}.$$

Construct
\begin{equation}
\label{eq:11.11}
\widehat{\Delta}^{J,L}_I:\widehat{Y}^I\to\widehat{Y}^J\wedge\widehat{Y}^L
\end{equation}
as follows. Let $$W^{J,L}_I$$ denote the smash product
$$\bigwedge_{\ell(J)+\ell(L)} W_i,$$
where $$W_i=\left\{\begin{array}{lll}Y_i&{\rm if}&i\in I-(J\cap
L)\\Y_i\wedge Y_i &{\rm if}&i\in J\cap L.\end{array}\right.$$

Note that if $J\cap L= \varnothing$, then $W_I^{J,L}=\widehat{Y}^I$.

Define $$\psi: \widehat{Y}^I\to W^{J,L}_I\quad{\rm as}\quad \psi=
\bigwedge_{i\in I}\psi_i$$ and $\psi_j:Y_j\to W_j$ as
$$ \psi_j=\left\{
\begin{array}{l}
\mbox{Id}\qquad{\rm if}\quad i\in I-(J\cap L)\\
\Delta_i:Y_i\to Y_i\wedge Y_i\quad{\rm if}\quad i\in J\cap L
\end{array}\right.
$$ where $\Delta_i$ is the reduced diagonal of $Y_i$.

Observe that the smash products $W^{J,L}_I$, and
$\widehat{Y}^J\wedge\widehat{Y}^L$ have the same factors, but in a
different order arising from the natural shuffles. So let
$$s: \widehat{Y}^J\wedge\widehat{Y}^L \to W^{J,L}_I$$
denote the natural homeomorphism given by a shuffle. Let
$$\theta:W^{J,L}_I \to \widehat{Y}^J\wedge\widehat{Y}^L $$
denote the inverse of $s$. Then define $$\widehat{\Delta}^{J,L}_I:
\widehat{Y}^I \to \widehat{Y}^J\wedge\widehat{Y}^L$$ as the
composition
\begin{equation}
\label{eq:new.28}
\widehat{Y}^I\stackrel{\psi}{\longrightarrow}W^{J,L}_I
\stackrel{\Theta}{\longrightarrow}\widehat{Y}^J\wedge\widehat{Y}^L.
\end{equation}

Further, observe that there is a commutative diagram
\begin{equation}
\label{eq:2.7} \xymatrix{ Y^{[m]} \ar[r]^{\Delta_{[m]}}
\ar[d]^{\widehat{\Pi}_J} &Y^{[m]} \wedge Y^{[m]}
\ar[d]^{\widehat{\Pi}_J\wedge\widehat{\Pi}_L}\\
\widehat{Y}^I\ar[r]^{\kern-.7cm\widehat{\Delta}^{J,L}_I}&
\widehat{Y}^J\wedge\widehat{Y}^L. }
\end{equation}

Next specialize Definition \ref{defin:star.product} to the case of
$(\underline{X},\underline{A}) =
(\underline{X},\underline{X})=\{(X_i,X_i)\}^m_{i=1}$ with $K
=\varnothing$. In this case, let $$\mathcal
H^*(K;(\underline{X},\underline{A})) = \mathcal H^*(X^{[m]}).$$ Then
the $*$-product in $\mathcal H^*(X^{[m]})$ is given by the
composition
$$H^p(X^J)\otimes H^q(X^L)\stackrel{{\rho^*}}{\to}H^{p+q}(X^J\wedge X^L)
\stackrel{(\widehat{\Delta}^{J,L}_I)^*}{\to}H^{p+q}(\widehat{X}^I)$$
where $I=J\cup L$. The notation $u * v$ denotes the $*$-product of
classes $u\in H^p(\widehat{X}^J), v\in H^q(\widehat{X}^L)$; the
class $u * v$ is a product of these two cohomology classes.
Furthermore, $$\mathcal H^n(Y^{[m]})\cong \bigoplus_I
H^n(\widehat{Y}^I)$$ as given in Definition
\ref{defin:star.product}.

Diagram (\ref{eq:2.7}) implies that
\begin{equation}
\Pi^*_I(u *v)=(\Pi^*_Ju)\smallsmile(\Pi^*_Lv)
\end{equation}
where $\smallsmile$ is the cup-product in $H^*(X^{[m]})$. The
well-known splitting of the suspension of a product stated in
\cite{bbcg} gives an additive isomorphism
$$\eta:\mathcal H^*(X^{[m]})\to H^*(X^{[m]}).$$ Since ${\eta}|H^*(\widehat{YX}^I)=(\Pi_I)^*$, diagram (\ref{eq:2.7})
implies that $\eta$ is a ring isomorphism.

The next result, a special case of Theorem \ref{thm:1.4}, follows.
\begin{thm}\label{thm:ring.isomorphism}
Let $(\underline{X},\underline{X})=\{(X_i,X_i)\}^m_{i=1}$ with $K =
\varnothing$. Then, the mapping
$$\eta:\mathcal H^*(X^{[m]})\to H^*(X^{[m]})$$
is a ring isomorphism.
\end{thm}

This special case of Theorem \ref{thm:1.4} will be extended to pairs
$(\underline{X},\underline{A})=\{(X_i,A_i)\}^m_{i=1}$ in the next
section.

\section{\bf The partial diagonal in the smash moment-angle complexes}\label{sec:The partial diagonal in smash moment angle}

In this section, the partial diagonal $\widehat{\Delta}^{J,L}_I:
\widehat{Y}^I \to \widehat{Y}^J\wedge\widehat{Y}^L$ of section
\ref{sec:The partial diagonal in product spaces} is extended to
smash moment-angle complexes as follows. Let $K$ be a simplicial
complex with $m$ vertices,  $\sigma\in K$ simplices of $K$, and
$(\underline{X},\underline{A})=\{(X_i,A_i)\}$. The simplices of
$K_I$ are $\sigma\cap I$, and $(\underline{X},\underline{A})_I$ the
associated family. Using the notation of section \ref{Introduction},
$$\widehat{Z}(K_I;\,\,(\underline{X},\underline{A})_I)=
\bigcup\limits_{\sigma\in K}\widehat{Y}^I(\sigma\cap I) \subset
\bigcup\limits_{\sigma\in K}\widehat{X}^I(\sigma\cap I).$$ Similarly
$$Z(K;(\underline{X},\underline{A}))=\bigcup\limits_{\sigma\in
K}D(\sigma)=\bigcup\limits_{\sigma\in K}Y(\sigma)\subset
\bigcup\limits_{\sigma\in K}X^I(\sigma) .$$

Observe that to give maps out of
$$Z(K;(\underline{X},\underline{A}))=\bigcup\limits_{\sigma\in
K}D(\sigma) = \operatorname{colim} (D(\sigma)),$$ it suffices to
give compatible maps out of each space $D(\sigma)$. In addition, to
give maps out of
$$\widehat{Z}(K_I;\,\,(\underline{X},\underline{A})_I)=
\bigcup\limits_{\sigma\in K}\widehat{Y}^I(\sigma\cap I) =
\operatorname{colim} (\widehat{Y}^I(\sigma\cap I)),$$ it suffices to
give compatible maps out of each space $\widehat{Y}^I(\sigma\cap
I)$.  Observe that the maps
$$\widehat{\Delta}^{J,L}_I: \widehat{X}^I \to
\widehat{X}^J\wedge\widehat{X}^L$$ as given in the composition
(\ref{eq:new.28}) restrict to maps
\begin{equation}
\label{eq:partial.diagonal.for.Z}
\widehat{\Delta}^{J,L}_I:\widehat{Z}(K_I) \to \widehat{Z}(K_J)\wedge
\widehat{Z}(K_L).
\end{equation}

\section{\bf Proof of Theorem \ref{thm:1.4}}\label{sec:Proof of Theorem (1.4)}

The proof of Theorem \ref{thm:bbcg} in \cite{bbcg} gives that
suspending and adding the maps $\widehat{\Pi}_I$ provides a map
$$\Sigma(\vee\widehat{\Pi}_I):\Sigma Z(K)\to \Sigma\vee_{I\subseteq [m]}\widehat{Z}(K_I)$$
which is a homotopy equivalence. Furthermore, each map
$$\widehat{\Pi}_I: Z(K)\to \widehat{Z}(K_I)$$ induces a morphism of
cohomology algebras while the sum of these maps induces
$$\eta:\mathcal H^*(K;(\underline{X},\underline{A}))\to H^*(Z(K;(\underline{X},\underline{A}))$$
which is an additive isomorphism. Since $\eta$ restricted to
$H^*(\widehat{K}_I)$ is ${\widehat{\Pi}_I}^*$, this implies that
$\eta$ is an algebra isomorphism, the statement of Theorem
\ref{thm:1.4}.

{\bf Note:} It is unnecessary to prove the associativity and the
graded commutativity of the $*$-product directly. Those properties
are a consequence of the isomorphism in Theorem \ref{thm:1.4}.

\section{\bf Proof of Theorem \ref{thm:1.5}}\label{sec:Proof of Theorem (1.5)}

If $Y_i$ is a suspension space then the reduced diagonal
$\Delta_i:Y_i\to Y_i\wedge Y_i$ is null-homotopic. Thus if $J\cap
L\neq \phi$, the map
$\widehat{\Delta}^{J,L}_I:\widehat{Z}(K_I)\to\widehat{Z}(K_J)\wedge\widehat{Z}(K_L)$
is null-homotopic. Theorem \ref{thm:1.5} follows.

\section{\bf Proof of Theorem \ref{thm:1.6}}\label{sec:Proof of Theorem (1.6)}

The first part of the Theorem is a consequence of the fact that
there is a homotopy equivalence
$$\widehat{Z}(K_I;(\underline{\Sigma^TX};\underline{\Sigma^TA})_I)\to
\Sigma^{d(T_I)}\widehat{Z}(K_I;(\underline{X},\underline{A})_I)$$
where for each $T=(t_1,\ldots,t_m)$, $$d(T_I)=\Sigma_{u\in I}t_i,$$ so
$$\widehat{Z}(K_I;(\underline{\Sigma^TX},\underline{\Sigma^TA})_I)$$
and
$$\widehat{Z}(K_I;(\underline{\Sigma^{T'}X},\underline{\Sigma^{T'}A}))$$
are stably equivalent.

The second part follows from the fact that all possible signs that
may come from permutations will have the same signs as $T\equiv T'$
mod $2$.

\section{{\bf The partial diagonal for $(CX,X)$} }\label{sec:CX.X}

The purpose of this section is to prove Theorem
\ref{thm:cup.products.for.CX.X} which exhibits the algebra structure
for the cohomology of $Z(K;(\underline{CX},\underline{X})).$ Let $K$
be an abstract simplicial complex with $m$ vertices. Given families
of based CW-pairs
$(\underline{Y},\underline{B})=\{(Y_i,B_i)\}^m_{i=1}$, and
$(\underline{W},\underline{C})=\{(W_i,C_i)\}^m_{i=1}$, consider the
natural shuffle map

\begin{equation}\label{eq:product.shuffle}
    \begin{array}{l}
    \mbox{shuff}:(Y_1\times W_1)\times \cdots \times (Y_m\times  W_m)\to
(Y_1\times \cdots \times Y_m) \times (W_1\times  \cdots\times W_m)
        \end{array}
\end{equation} restricted to

\begin{equation}
\label{eq: restrict}
\begin{array}{l}
\mbox{shuff}:Z(K;(\underline{Y\times W},\underline{B\times C})) \to
Z(K;(\underline{Y},\underline{B})) \times
Z(K;(\underline{W},\underline{C})).\\
\end{array}
\end{equation}

Let $(Y,B) = (D^1, S^0)$, and
$(\underline{W},\underline{C})=(\underline{X},\underline{X})$. There
is an induced map
\begin{equation}
\label{eq: distributivity over sums and iterations}
    \begin{array}{l}
\gamma:\widehat{Z}(K;(\underline{D^1\wedge X},\underline{S^0\wedge
X})) \to \widehat{Z}(K;(\underline{D^1},\underline{S^0})) \wedge
\widehat{Z}(K;(\underline{X},\underline{X}))\\
\end{array}
\end{equation} where $Z(K;(\underline{X},\underline{X}))= X_1\wedge \cdots \wedge X_m =
\widehat{X}^{[m]}$. Furthermore,
$$\widehat{Z}(K;(\underline{D^1\wedge X},\underline{S^0\wedge X}))=
\widehat{Z}(K;(\underline{CX},\underline{X})).$$ The next Lemma
follows by inspection.

\begin{lem}\label{lem:decomposition} Let $K$ be an abstract simplicial complex with $m$ vertices.
The natural map induced by the shuffle
$$\gamma:\widehat{Z}(K;(\underline{D^1 \wedge X},\underline{S^0 \wedge X})) \to
Z(K;(\underline{D^1},\underline{S^0})) \wedge
Z(K;(\underline{X},\underline{ X}))$$ is a homeomorphism. Thus
$\widehat{Z}(K;(\underline{D^1\wedge X},\underline{S^0\wedge X}))=
\widehat{Z}(K;(\underline{CX},\underline{X}))$ is naturally
homeomorphic to
$$\widehat{Z}(K;(\underline{D^1},\underline{S^0}))\wedge
\widehat{X}^{[m]}.$$
\end{lem}

Let $\widehat{Z}(K_I)$ denote
$\widehat{Z}(K_I;(\underline{CX},\underline{X}))$ and let
$$\gamma_I: \widehat{Z}(K_I;(\underline{CX},\underline{X})) \to
\widehat{Z}(K_I;(\underline{D^1},\underline{S^0})) \wedge
\widehat{X}^{K_I}$$ denote the homeomorphism of Lemma
\ref{lem:decomposition}. Observe that if $ J \cup L = I$, there is a
homotopy commutative diagram
\[
\begin{CD}
 \widehat{Z}(K_I) @>{\ \  \  \   \   \widehat{\Delta}^{J,L}_I} \  \  \  \  >> \widehat{Z}(K_J) \wedge \widehat{Z}(K_L)  \\
@V{\gamma_I}VV          @VV{\gamma_J \wedge \gamma_L}V          \\
\widehat{Z}(K_I;(D^1,S^0))\wedge\widehat{X}^{I}
@>{\widehat{\Psi}^{J,L}_I(\underline{CX}, \underline{X})}>>
\widehat{Z}(K_J;(D^1,S^0))\wedge\widehat{X}^{J}\wedge
\widehat{Z}(K_L;(D^1,S^0))\wedge\widehat{X}^{L}
\end{CD}
\] where the map $\widehat{\Delta}^{J,L}_I:\widehat{Z}(K_I) \to \widehat{Z}(K_J) \wedge
\widehat{Z}(K_L)$ is given by (\ref{eq:partial.diagonal.for.Z}).

Specialize to $(CX,X)= (D^1,S^0)$ with  $\widehat{Z}(K_I)=
\widehat{Z}(K_I;(\underline{CX},\underline{X}))$ as well as to the
unique pointed homeomorphism $X\wedge \cdots \wedge X \to X = S^0$
to obtain a homotopy commutative diagram

\[
\begin{CD}
 \widehat{Z}(K_I) @>{\  \   \    \   \    \widehat{\Delta}^{J,L}_I \   \   \   }>> \widehat{Z}(K_J) \wedge \widehat{Z}(K_L)  \\
@V{\alpha_I}VV          @VV{\alpha_J \wedge \alpha_L}V          \\
\widehat{Z}(K_I;(D^1,S^0))@>{\widehat{\Psi}^{J,L}_I(D^1,S^0)}>>
\widehat{Z}(K_J;(D^1,S^0)) \wedge \widehat{Z}(K_L;(D^1,S^0)).
\end{CD}
\] This diagram is used to address the general case.

Consider $(\underline{CX},\underline{X})$ together with the map
\[
\begin{CD}
\widehat{Z}(K_I;(D^1,S^0))\wedge\widehat{X}^{I}
@>{\widehat{\Psi}^{J,L}_I(\underline{CX},\underline{X})}>>
\widehat{Z}(K_J;(D^1,S^0))\wedge\widehat{X}^{J}\wedge
\widehat{Z}(K_L;(D^1, S^0))\wedge\widehat{X}^{L}.
\end{CD}
\] Observe that the map $\widehat{\Psi}^{J,L}_I(\underline{CX},\underline{X})$ is given by the composite
\[
\begin{CD}
\widehat{Z}(K_I;(D^1,S^0))\wedge\widehat{X}^I
@>{\widehat{\Psi}^{J,L}_I(D^1,S^0)\wedge \widehat{\Delta}^{J,L}_I
}>> \widehat{Z}(K_J;(D^1,S^0)) \wedge
 \widehat{X}^J \wedge \widehat{Z}(K_L;(D^1,S^0)) \wedge \widehat{X}^L
\end{CD}
\]

\noindent with

\[
\begin{CD}
\widehat{Z}(K_J;(D^1,S^0)) \wedge \widehat{X}^J \wedge
\widehat{Z}(K_L;(D^1,S^0)) \wedge \widehat{X}^L @>{ 1 \wedge \tau
\wedge 1}>> \Sigma|K_J|\wedge \widehat{X}^J \wedge \Sigma|K_L|\wedge
\widehat{X}^L
\end{CD}
\]

\noindent where

$$\tau: \widehat{Z}(K_L;(D^1,S^0))
 \wedge\widehat{X}^J \to \widehat{X}^J \wedge \widehat{Z}(K_L;(D^1,S^0))$$ is the
natural map which swaps factors. This last assertion follows by
inspection of the definitions. This suffices to prove Theorems
\ref{thm:cup.products.for.CX.X}, and
\ref{thm:E.cup.products.for.CX.X}.

\section{{\bf Acknowledgements} }\label{Acknowledgments}. 

The authors are grateful to Taras
Panov, Graham Denham, Matthias Franz, and Peter Landweber for their
interesting questions and helpful suggestions. This paper depends
heavily on their questions, and is better because of their
suggestions. The first author was supported in part by the award of
a Rider University Research Fellowship and the third author was
partially supported by DARPA grant number 2006-06918-01.

\bibliographystyle{amsalpha}

\end{document}